\documentclass[11pt]{amsart}
\usepackage{epsfig,amsmath,amssymb,amsthm,epsf}
\usepackage{graphicx,array}

\DeclareMathOperator{\Prob}{P}

\DeclareMathOperator{\E}{E}

\newtheorem{proposition}{Proposition}

\newtheorem{theorem}{Theorem}
\newtheorem{claim}{Claim}

\begin{document}
\title{{The Escape model on a homogeneous tree.}}
\author{{George Kordzakhia}}
\address{University of California\\
Department of Statistics\\
Berkeley CA }
\email{kordzakh@stat.berkeley.edu}
\maketitle
\begin{abstract}
There are two types of particles  interacting on a homogeneous tree of degree $d+1$. The particles of the first type  colonize  the empty space with exponential rate $1$, but cannot take over the vertices that are occupied  by the second type.    
The  particles of the second type  spread with exponential rate $\lambda$.
They colonize   the neighboring vertices that are either vacant or occupied by the  representatives of the opposite type,  and  annihilate  the particles of the type $1$ as they reach them. There exists  a critical value $\lambda_c=(2d-1)+\sqrt{(2d-1)^2-1}$ such that the first type survives with positive probability for $\lambda < \lambda_c$, and dies out
with probability one for $\lambda > \lambda_c$.
We also find the growth profile which characterizes  the rate of growth  of the type $1$ in the space-time on the event of survival.
\end{abstract}
\section{Introduction}
We consider a model of a predator-prey type which we call the Escape  model.
There are two entities  growing on the vertices of a homogeneous tree
$\mathcal{T}_d=\mathcal{T}$ of degree $d+1$.
The entities may be thought of as  biological species, political parties or
manufacturers competing on a market.  
The second  entity dominates the first in the sense
that the representatives of the second entity can take over
the vertices occupied by the representatives of the first entity
but not vice versa. It is also  assumed that the second entity  grows faster. 
We are interested in the possibility of the long-term coexistence of 
the species which occurs  when the first species survives.

At each time $t$  each site of  the tree is 
occupied by at most one  representative  of either the two entities. 
 We  refer to  the representatives of the  entities as 
 particles of  types $1$ and $2$ respectively.
If a site is not occupied, it is said  to be {\it vacant}. 
 The dynamics of the process 
is specified by exponential transition rates.
A vacant site gets  colonized  by a particle of type  $1$
 with exponential rate equal to the number of nearest neighbors of type $1$.  
The sites  that are  either vacant or occupied by type $1$  flip to $2$
 with the  rate   $\lambda > 1$ times  the number of neighbors of type $2$. 
If a vertex is  occupied by a particle of type $2$, the particle stays 
at the vertex forever. We assume that at time zero  there are finitely many particles 
of each  type. 

The sets of sites occupied by particles of types $1$ and $2$ 
at time $t \ge 0$ are denoted by $A(t)$ and $B(t)$ respectively.
The event $\cap_{t \ge 0} \{ A(t) \not= \emptyset\}$ is referred 
to as the event of {\it survival} of type $1$. 
If at time zero the cluster $A(0)$  is surrounded by $B(0)$, 
then all particles of type $1$  eventually die out with probability one.
Thus, we are only interested in  the  initial configurations for which   
there exist a vertex $x$ and an infinite  geodesic segment $\gamma_{x,\infty}$
 such that $x \in A(0)$ and $\gamma_{x,\infty} \cap B(0)= \emptyset$.
All such configurations are referred to as {\it  non-trivial} configurations.
The question of interest is for what values of $\lambda$ type $1$  survives 
with positive probability in the long term. 
\begin{theorem} \label{Escape:MainTheorem}
There exists  a critical value  $ \lambda_c=(2d-1)+\sqrt{(2d-1)^2-1}$ such that,  
for all $\lambda  \in   (1, \lambda_{c})$ and for all
non-trivial  finite configurations $(A(0), B(0))$,  type $1$ survives with positive 
probability. 
For all $\lambda  \in  (\lambda_{c}, \infty)$, type $1$ dies out with probability one.
\end{theorem} 
For $c>0$ denote by $M_n(n/c)$ the number of vertices $x \in  A(n/c)$ at distance $n$ 
from the  root $\rho$. From the results of  Sections  \ref{Two independent Richardson models.} and \ref{section:Escape}, it follows that 
 \[ \lim \frac{1}{n}\log\left(\E M_{n}(n/c)\right)= -g(c). \] 
where 
\[
g(c)=\left\{ 
\begin{array}
{ l@{\quad: \quad} l}
\left( \lambda/c-\log (\lambda/c)-1\right)-\log d &  0<c \le 1 \\
\left(\lambda/c-\log (\lambda/c)-1\right)+ \left(1/c-\log (1/c)-1 \right)-\log d &  1 < c < \lambda \\
 \left(1/c-\log(1/c)-1 \right)-\log d  &  c \ge  \lambda
\end{array}
\right.
\]
The function $g(c)$ is referred to as  the {\it growth profile} of type $1$.  
The growth profile was introduced in Lalley \cite{lalley2} (in a slightly different form)  
to study  the weakly supercritical  contact process on a homogeneous tree.
The function $g(c)$ has a unique minimum at $c_0$, is strictly decreasing on $(0, c_0)$ and strictly increasing on $(c_0, \infty)$.
For  all  $\lambda  \in   (1, \lambda_{c})$, we have  $g(c_0)<0$, and
let $r_1$ and $r_2$, with $0<r_1<r_2$,  be the two solutions of $g(c)=0$.
\begin{theorem} \label{Escape:NumberOccupied}
Let  $\lambda  \in   (1, \lambda_{c})$. 
For every $\epsilon>0$ and all large $t$,  the particles of type $1$ are concentrated in
the annulus  of radii $(r_1- \epsilon)t$ and $(r_2+\epsilon)t$ centered at the root. 
For every  $c \in (r_1, r_2)$, almost surely on the event of survival of type $1$,
\[ \lim \frac{1}{n}\log\left(M_{n}(n/c)\right)= -g(c). \] 
\end{theorem}
To investigate the Escape model, we consider  the simplest growth model, 
the Richardson model, presented in Section \ref{section:Richardson}.                     
\section{Preliminaries.}
\subsection{ A homogeneous tree}
\indent  A homogeneous tree   $\mathcal{T}_d=\mathcal{T}$ of degree $d+1$ is an infinite tree
such that every vertex has exactly $d+1$ nearest neighbors.
A distinguished vertex is called {\it root} and denoted by $\rho$. 
For every two vertices $x$ and $y$ of the tree, denote  by $|x,y|$ 
the number  of edges in the shortest  path  from $x$ to $y$ (the path without loops). 
If $x$ is the root, then we  simply write $|y|$. Note that $| \cdot , \cdot |$ is a metric
on $\mathcal{T}$. \\
For  every vertex $x \in \mathcal{T}$, a {\em geodesic segment} $\gamma_{x,\infty}$ is an infinite path in $\mathcal{T}$ beginning at $x$ and having  no loops.
Define by $\mathcal{T}_{+}(x)$, the set of all vertices $y$ such that the shortest path connecting $y$  with the root $\rho$ goes through $x$.
Consider  a geodesic segment $\gamma_{x,\infty}$  such that $\gamma_{x,\infty} \in \mathcal{T}_{+}(x) $.  Given an integer $m$, consider also  a  sequence of vertices $ (y_k)_{k \ge 0}$ on the geodesic segment 
$\gamma_{x,\infty}$  such that $|y_k|=mk+|x|$  for all integers $k \ge 0$. For every  $k \ge 1$,  $y_{k-1}$  is called the { \it $m$-predecessor} of $y_k$ in $\mathcal{T}$. \\ 
 Let  $D(x,r)=\left\{ y\in \mathcal{T}:\ |x,y|\le r \right\}$  be the
 closed  disk of radius $r \in (0,\infty)$  centered at $x$, and let $C(x,r)=\left\{ y\in
   \mathcal{T}:\ |x,y|= r \right\}$ be  the circumference of that
 disk. If $x$ is the root, then we write $D_r$ and $C_r$.
Note that for all integers $n \ge 1$, the number of vertices in $C_n$ is $(d+1)d^{n}$.
\subsection{Construction of the Escape process.}
\indent The Escape process can be built using a percolation structure as follows.
For each ordered pair  of neighboring vertices $x$ and $y$ in 
$\mathcal{T}$,  define two  independent Poisson processes with respective rates
 $1$ and $\lambda-1$,  and  respective  occurrence times $\left\{ T_n^{x,y}: n \ge 1  \right\}$
and $\left\{ U_n^{x,y}: n \ge 1 \right\}$.  Set $T_0^{x,y}=0$ and $U_0^{x,y}=0$,
and make these Poisson processes independent from pair to pair.
Consider  $\mathcal{T} \times \mathbb{R}_{+}$.  Arrows  are drawn from 
$x$ to $y$  at the occurrence times  $T_n^{x,y}$
and  $U_n^{x,y}$.
We say that there is a {\it directed path} in  $\mathcal{T} \times \mathbb{R}_{+}$ 
from $(x_0,s_0)$ to $(x_{n},s_{n+1})$ if there is a sequence 
of times $s_0 <  s_1<..<s_{n+1}$  and sequence of vertices $x_0,x_1,..,x_n$ so that
for each $j,\ 1<j \le n$, there is an arrow  from $x_{j-1}$ to $x_j$ at time $s_j$.
A {\it type $1$ path} is a directed  path that uses only  arrows
generated by the Poisson processes $T_n$. 
Fix an initial configuration $A(0), B(0)$ and erase all arrows that lie only on paths that begin
at points $(x,0)$ such that $x \not \in  (A(0) \cup B(0))$. 
For every vertex $y$ we say that  $y \in B(t)$ if and only if there is  a vertex $z \in B(0)$ and  
directed path from  $(z,0)$ to $(y,t)$ in the modified percolation structure.
Define $A(t)$ to be the set of vertices $y$
such that $y \not \in B(t)$ and  there is a type $1$ path (in the modified percolation structure) that ends at $(y,t)$ and starts at $(z,0)$ for some  $z \in A(0)$.
\section{ The Richardson model on a homogeneous tree} \label{section:Richardson}
\indent The Richardson process on $\mathcal{T}$  with 
parameter $\lambda>0$  is a continuous time  Markov process $R(t)$ on
the set of finite subsets of $\mathcal{T}$. We say that a vertex  $x$
is { \it infected} (or occupied) at time $t$ if $x \in R(t)$, and  is vacant otherwise.
The process  develops according to the following rules: if a  vertex gets infected, it stays infected forever and starts infecting unoccupied neighboring vertices with rate $\lambda$, i.e. that the infection times have exponential distributions with parameter $\lambda$. The infection times are all independent. Consequently,  a vacant site becomes  infected with  the  rate  
\[\lambda \cdot \mbox{(number of infected neighbors)}. \]
The model with parameter $\lambda=c$  can be obtained from the model with parameter $\lambda=1$  by running the later  process at speed $c$. Therefore we treat just the case $\lambda=1$.  \\ 
\indent Without loss  of generality consider the initial configuration where at time zero the only
 occupied site is the root $R(0)=\{\rho\}$.
The  main questions  were how fast the infected set grows  and what limiting shape the infected  set has.  For the tree $\mathcal{T}$ the number of vertices in the disk
$D_n$ grows exponentially with $n$ (whereas for $\mathbb{Z}^d$  it has
polynomial growth). Consequently, the behavior of the Richardson process on
$\mathcal{T}$ is different from the behavior  of the model on 
the integer lattice (see \cite{richard} and \cite{cox}). 
The infected region $R(t)$ still  grows linearly with time, but
there are constants  $a$ and $b$, with $a<b$, such that as $t$ goes to
infinity we can classify two subregions: a completely infected
subregion, having approximately the  shape of a ball $D_{ta}$, and a
partially infected subregion, having approximately  the shape of a
ring $D_{tb} \setminus D_{ta}$. Consequently, there are two  speeds:
the speed of invasion $b$, indicating how fast the infection spreads,
and the speed of occupation $a$,  governing the rate of growth of the
region that is completely covered by the infection.   
\begin{proposition} \label{Rich:TreeShape}
Let $d \ge 2$ be an integer. Let 
\begin{equation} \label{Rich:TreeShape:function:f}
f(c)=\frac{1}{c}-\log \frac{1}{c}-1 -\log d
\end{equation} 
for $c\in (0,\infty)$,  and~let~$a$~and~$b$, with $0<a<1<b$, be the two roots 
of the equation 
\[ f(c)=0. \]  
Then $a=\sup\{a'\}$ and $\ b=\inf\{b'\}$, where the  $\sup$ and $\inf$ are taken over all $a'$ and $b'$  satisfying 
\[\Prob[\exists \ \mbox{a random} \ \tau < \infty \ \mbox{such that},\ \forall\  t> \tau ,\ \ D_{ta'} \subset R(t) \subset D_{tb'}]=1. \] 
$ \mbox{ As } \  d \rightarrow  \infty,\ \mbox{we have} \ a(d) \log d \rightarrow 1 \ \mbox{and}\ b(d)/d \rightarrow e.$ \\
\end{proposition}
Let $N_n(t)$ be the number of vertices at  distance $n$ from the root that  are  infected at time $t$, and let  $F_n(t)$ be the number of  vertices in $C_n$ that are {\em not} infected at time $t$. 
We  compute asymptotic values of  $N_n(t)$ and  $F_n(t)$  as $n$ goes to infinity and $t=n/c$ for different values of  $c$. 
\begin{proposition} \label{Rich:NumberOccupied} For all $ c\in (1, b)$, 
\begin{equation}  \label{Eq:NumberOccupied}
\lim\frac{1}{n}\log \left( N_{n}({n/c}) \right)=-f(c)>0 \ \mbox{a.s.}. 
\end{equation}
For all $c\in (a, 1)$,
\begin{equation} \label{Eq:NumberVacant}
 \lim\frac{1}{n}\log \left( F_{n}({n/c})\right)=-f(c)>0 \ \mbox{a.s.}.  
\end{equation}
\end{proposition}
\begin{proof} [Proof of Proposition \ref{Rich:TreeShape}] 
It is enough to verify that for every   $\epsilon \in (0,a)$ \\
\begin{equation}  \label{Rich:TreeShape1}
 \Prob[\exists   \tau<\infty  \ \mbox{such that} \  \forall  t>\tau, \ D_{t(a-\epsilon)} \subset R(t) \subset D_{t(b+\epsilon)}]=1, 
 \end{equation}
\begin{equation} \label{Rich:TreeShape2}
\Prob[\exists   \tau<\infty \ \mbox{such that} \ \forall t>\tau,\  R(t) \subset D_{t(b-\epsilon)}]=0, 
\end{equation}
\begin{equation} \label{Rich:TreeShape3}
\Prob[\exists   \tau<\infty \ \mbox{such that} \ \forall  t>\tau, \ D_{t(a+\epsilon)} \subset R(t)] =0.
\end{equation} \\
Proof of equation (\ref{Rich:TreeShape1}).  For every vertex $x$ on the tree, 
let $T(x)$ be the (random) time at which the vertex gets infected. 
Consider an arbitrary vertex at distance $n$ from the root, and denote it  by $x_n$.
Fix $c>1$. Then we have the following estimate: \\
\begin{gather*}
  \Prob\left[x_n \in R(n/c)\right] = \Prob\left[T(x_n) \le \frac{n}{c}\right]=\exp\left\{-n\left(\frac{1}{c}-\log\frac{1}{c}-1\right)+o(n)\right\} \\
 \intertext{where}
o(n)/ n \rightarrow 0 \ \ \mbox{as} \ \  n \rightarrow \infty.
\end{gather*}
The estimate  follows from  Cram\'{e}r's theorem for i.i.d. random variables (see \cite{dembo}) and from the fact that $T(x_n)$ is distributed as a sum of $n$ i.i.d. exponentials with parameter $1$. 
Observe that
\[
 \E N_{n}(n/c)=(\mbox{the total number of vertices in }C_n)\cdot \Prob(x_n \in R(n/c))  \]
\[  =\frac{d+1}{d}d^n\exp\left\{-n\left(\frac{1}{c}-\log \frac{1}{c}-1 \right)+o(n)\right\} \] 
\begin{equation} \label{Rich:ExpOccupied}
=\exp\{-nf(c)+o(n)\}
\end{equation} \\ 
where $f$  was defined in (\ref{Rich:TreeShape:function:f}). \\
\indent  Similarly, for $c<1$ we have:
\begin{equation*}
 \Prob[x_n \notin R(n/c)]= \Prob\left[T(x_n) \ge \frac{n}{c}\right] 
=\exp\left\{-n\left(\frac{1}{c}-\log \frac{1}{c}-1\right)+o(n)\right\}. 
\end{equation*} 
Hence
\begin{equation}  \label{Rich:ExpVacant}
\E F_{n}(n/c)=\exp\{-nf(c)+o(n)\}. 
\end{equation}
Observe that $f$ is 
strictly decreasing on $(0,1)$ and strictly increasing on $(1,\infty)$, 
with   unique minimum at $c=1$.  Moreover $f(1)<0$ and $f(0+)=f(\infty-)=\infty$. 
Thus, there are just two roots  $a$ and $b$  of the equation $f(c)=0$,  
such that $a<1<b$. Hence, for every $\epsilon >0$, we have that \\
\[ \Prob\left[ \exists x \in C_n: \  T(x) \le \frac{n}{b+\epsilon}\right] \le \E N_{n}(n/(b+\epsilon)). \] 
\noindent Since $f(b+\epsilon)>0$, 
by (\ref{Rich:ExpOccupied}) the upper bound decays exponentially with $n$. 
Consequently  the left  side is summable, and  by Borel-Cantelli lemma  we have  \\
\[ \Prob\left[\exists  N < \infty: \ \forall  n > N \ \mbox{and} \ \forall  x \in C_n,\ T(x)>\frac{n}{b+\epsilon}\right]=1. \] 
\noindent Finally, observe that this is equivalent to
\begin{equation} \label{Rich:TreeShape1Upper}
 \Prob\left[ \exists   \tau < \infty:  \  \forall
   t>\tau,\  R(t) \subset D_{t(b+\epsilon)} \right]=1
\end{equation}
 (the  events are identical). \\
\indent Analogously,  for all $\epsilon \in(0,a)$,
\begin{equation} \label{Rich:TreeShape1Lower}
 \Prob\left[\exists \tau < \infty: \  \forall   t>\tau ,\ D_{t(a-\epsilon)} \subset R(t) \right]=1. \end{equation} 
 To  prove (\ref{Rich:TreeShape1Lower}),  note that  
\[ \Prob\left[ \exists x \in C_n:\  T(x) > \frac{n}{a-\epsilon} \right] \le \E F_{n}(n/(a-\epsilon)), \]
and  it decays exponentially  by (\ref{Rich:ExpVacant}). 
Applying  Borel-Cantelli lemma one more time, we get  
\[\Prob \left[\exists N< \infty : \ \forall  n > N  \ \mbox{and} \  \forall  x \in C_n, \ T(x) \le \frac{n}{a-\epsilon}\right]=1 \] 
which implies (\ref{Rich:TreeShape1Lower}).
Obviously, (\ref{Rich:TreeShape1Upper}) and (\ref{Rich:TreeShape1Lower}) together  are equivalent to (\ref{Rich:TreeShape1}).\\
 \indent Equations  (\ref{Rich:TreeShape2}) and (\ref{Rich:TreeShape3}) are direct consequences of  Proposition \ref{Rich:NumberOccupied}. As a heuristic argument, note that (\ref{Rich:ExpOccupied}) implies that, for every $c\in (1,b),  \  \E N_{n}(n/c)$ grows exponentially,  and, similarly, (\ref{Rich:ExpVacant}) implies that 
for every $c \in (a,1),\   \E F_{n}(n/c)$ grows exponentially.  
These observations suggest that (\ref{Rich:TreeShape2}) and (\ref{Rich:TreeShape3}) should be true. \\
\indent The limits 
\[ \lim_{d \rightarrow \infty}  a(d) \log d =1 , \]
\[ 
 \lim_{d \rightarrow \infty}  b(d)/d = e
\]
immediately follow  from the fact that $a$ and $b$ are the roots of 
\[ \frac{1}{c}-\log\frac{1}{c} - 1 -\log d = 0. \] 
The phenomenon  is easily anticipated.
Since for larger $d$'s there are more directions for the infection to spread around, it is natural that the invasion speed  is strictly increasing with $d$.  For the same reason, the occupation speed decreases to zero
(the number of vertices in $D_n$  grows unboundedly with $d$).
\end{proof} 
\begin{proof} [Proof of Proposition \ref{Rich:NumberOccupied}]
We only prove (\ref{Eq:NumberOccupied}). The proof of (\ref{Eq:NumberVacant}) is identical.
First we claim that for any $\epsilon >0$ 
\[ \limsup \frac{1}{n}\log\left(N_{n}(n/c)\right) \le -f(c)+\epsilon \] \ 
\noindent almost surely. 
By Markov's  inequality and (\ref{Rich:ExpOccupied}),
\begin{equation*}
\Prob\left[\frac{1}{n} \log \left(N_{n}(n/c)\right) > -f(c)+\epsilon \right] 
= \Prob\left[ N_{n}(n/c) > \exp\left\{ n(-f(c)+\epsilon) \right\} \right] \le 
\end{equation*}
\begin{equation*}
\le \E\left[ N_{n}(n/c) \right] \cdot \exp\left\{ nf(c)-n\epsilon \right\}  = \exp\left\{ -n\epsilon + o(n) \right\},
\end{equation*}
and the claim follows by Borel-Cantelli lemma.\\
\indent Let $\alpha=\frac{1}{c}$. 
To finish the proof, 
it is enough to show that, for an arbitrarily small $\epsilon>0$, with probability $1$
\begin{equation} \label{Rich:NumberOccupiedLower}
\liminf \frac{1}{n}\log\left(N_{n}(\alpha n)\right) \ge -f\left( 1/\alpha \right)-\epsilon 
\end{equation} 
Recall that $1<c<b$, so $f \left( 1/ \alpha \right)<0$.  Assume that $\epsilon$ is small enough to satisfy
 $f \left( 1 / \alpha \right)+\epsilon < 0$. To make the notation less complicated, let 
\[ \mu=\mu(\epsilon)=\exp\left\{-\left(f\left(1 / \alpha \right)+\epsilon\right) \right\}. \] 
By the~continuity of $f$, there exists an $\epsilon_1 >0$ such that    
\[ f \left( 1 / (\alpha-2\epsilon_1) \right) <  f\left(1 / \alpha \right)+\epsilon < 0. \] 
Then,  for every fixed $w>1$, there exists an integer $m>0$  large enough such that
\[ \E N_m((\alpha-2\epsilon_1)m)
> w\cdot \exp\left\{ -m\left( f\left(1 /\alpha  \right)+\epsilon\right) \right\} \]
\begin{equation} \label{Rich:NumberOccupied:Constant>1}
 =w \cdot \mu^m > 1. 
\end{equation}
 Choose $m$ to satisfy  (\ref{Rich:NumberOccupied:Constant>1}). Fix an arbitrary vertex $x$ of the tree and consider a geodesic segment  $\gamma_{x,\infty} \in \mathcal{T}_{+}(x)$. 
Consider a  sequence of vertices $ ( y_i )_{i \ge 0}$ 
on the geodesic segment $\gamma_{x,\infty}$  such that $ |y_i|=mi+|x|$  for all integers $i\ge 0$. 
Note that infection times $T(y_i)$ are increasing in $i$.
For every  pair of non-negative  integers $n_1$ and $n_2$ such that $ n_1 < n_2 $,  
say that $y_{n_2}$ is an { \it $m$-descendant} of~$y_{n_1}$ if, for all integers  $i \in [n_1,\ n_2)$,  
\[ T(y_{i+1})-T(y_{i})<\ (\alpha-2\epsilon_1)m. \]  
Define 
 \[ \mathcal{Z}_{k}(x) =\left\{ z \in C_{km+|x|}: \  z \ \mbox {is an} \ m\mbox{-descendant of}\  x \right\}, \]
\[ Z_{k}{(x)} = \mbox{ cardinality of}\ \mathcal{Z}_k{(x)}. \] \\
 Note that $(Z_{k}{(x)})_{k\ge 0}$ is a Galton-Watson process with mean offspring number
 \begin{equation} \label{Rich:ExpGaltonWatson}
\E Z_{1}{(x)}  > w \cdot \mu^m > 1. 
\end{equation}
\begin{claim}
For every vertex $x$, almost surely on the event of  ${(Z_{k}{(x)})}_{k\ge 0}$ survival, 
there exists a (random)  $K <+\infty$  such that,  for all $k>K$,  
\begin{equation} \label{Rich:NumberOccupied:GaltonWatson} 
N_{mk+|x|} ((\alpha -\epsilon_1)(mk+|x|)) \ge Z_{k}{(x)} \ge \mu^{mk+|x|}\cdot d^m.  
\end{equation}   
\end{claim}
\begin{proof}
The first inequality follows from the fact that for sufficiently large $k$'s  
\[ (\alpha-\epsilon_1)(mk+|x|) \ge (\alpha-2\epsilon_1)mk+ T(x). \]
To obtain the second inequality,  observe that, for all large enough $k$,
 $w^k \ge d^m \mu^{|x|}$, and hence, $\E Z_k(x)> (\mu^m \cdot w)^k 
\ge \mu^{mk+|x|}\cdot d^m$.
A standard theorem from the elementary  theory of Galton-Watson processes
states that if $\E Z_{1}{(x)}  > 1$  and  the variance of $Z_{1}{(x)}$ is finite, 
then on the event of survival $Z_{k}^{(x)} / (\E Z_{1}^{(x)})^k$ converges almost surely 
to a positive random variable. Thus, the second inequality is obtained by 
direct application of (\ref{Rich:ExpGaltonWatson}).
\end{proof}
Therefore, on the event  of non-extinction of  ${(Z_{k} {(x)})}_{k\ge 0}$, (\ref{Rich:NumberOccupied:GaltonWatson})  is true and implies (\ref{Rich:NumberOccupiedLower})  for  $n$'s  of the form $mk+|x|$. 
To establish the result for all positive integers, consider $mk+|x| < n < m(k+1)+|x|$.
Observe that $K$ might be also chosen large enough  that, for all $k>K$, we have 
\[  (\alpha-\epsilon_1)(m(k+1)+|x|\ ) <  \alpha(mk+|x|\ ). \]  
By (\ref{Rich:NumberOccupied:GaltonWatson})  at time  $(\alpha-\epsilon_1)(m(k+1)+|x|)$  we have at least 
\[ \mu^{m(k+1)+|x|}\cdot d^m  \] infected vertices on level $m(k+1)+|x|$.
Since each particle can generate at most $d$ offspring, it follows that, for each integer $n$ satisfying  
$mk+|x| < n < m(k+1)+|x|$,  there are at least
\[ \mu^{m(k+1)+|x|\ }d^{n-(mk+|x|)} \] infected predecessors in $C_n$ at  time $(\alpha-\epsilon_1)(m(k+1)+|x|)$. Obviously,
\begin{align*}
 N_{n}(\alpha n) &> N_{n}(\alpha(mk+|x|)) \\
\rule{0pt}{20pt}
&> N_{n}((\alpha-\epsilon_1)(m(k+1)+|x| )) > \mu^{m(k+1)+|x|} > \mu^{n}. 
\end{align*}
Therefore we proved that, almost surely  on the event of the survival of $(Z_{k}{(x)})_{k \ge 0}$, (\ref{Rich:NumberOccupiedLower}) is true. 
 To show that (\ref{Rich:NumberOccupiedLower}) holds with probability one, observe that, for each  integer $j>0$, there are $(d+1)d^{mj-1}$  Galton-Watson processes $(Z_k(x))_{k \ge 0}$ with $|x|=mj$. Let $S_j$ be the event of non-extinction for at least one of the processes. Since the Galton-Watson processes are independent, the probability of  $S_j$  tends to 1 as $j$ tends to $+\infty$. Moreover, for each $j>0$,  $S_j \subset S_{j+1} $   which guarantees the almost sure  result.  
\end{proof}
\section{ Two independent Richardson models.} \label{Two independent Richardson models.}
To motivate the proof of Theorem \ref{Escape:MainTheorem}
consider two { \it independent} Richardson processes $R_1(t)$ and $R_{\lambda}(t)$ 
with respective rates $1$ and $\lambda>1$.  
The processes  have   initial configurations $R_1(0)=R_{\lambda}(0)=\{\rho\}$, and  
are built on a  homogeneous tree $\mathcal{T}$. \\
For $c \in (0,\infty)$,  estimate the expected number of vertices   
in  $(R_1(n/c) \setminus R_\lambda(n/c)) \cap C_n$, that is the number of vertices at
 distance $n$  from the root that are occupied by $R_1$ but not by $R_{\lambda}$ at 
time $n/c$.  Let  $x_n$ be a vertex with $|x_n|=n$, 
$u_n(n/c)=\Prob \left[ x_n \in \left\{R_1(n/c) \setminus R_\lambda(n/c)\right\} \right]$, and $V_n(n/c)= \#  \left\{ x_n: x_n \in R_1(n/c) \setminus R_\lambda(n/c) \right\}$.
 \\
\indent Case 1:  For  every $c \in (0,1]$,  
\[
u_n(n/c)=
\Prob \left[ x_n \not \in R_{\lambda}(n/c) \right]
\Prob \left[ x_n \in R_1(n/c)\right]=  
\]
\[
= \exp \left\{-n\left(\frac{\lambda}{c}-\log \frac{\lambda}{c}-1 \right) + o(n) \right\}.
\]
Recall that the number of vertices in $C_n$ is $\frac{d+1}{d}d^n$ and  define 
\[
g_1(c)=\left(\frac{\lambda}{c}-\log \frac{\lambda}{c}-1\right)-\log d.
\]
Then, for all $c \in(0,1]$,
 \[
\E  V_n(n/c) =
 \exp \left\{-n \cdot g_{1}(c)+ o(n) \right\}.
\]

\indent Case 2: For every  $c \in (1,\lambda)$
\[
u_n(n/c)=
\Prob \left[ x_n \not \in R_{\lambda}(n/c) \right]
\Prob \left[ x_n \in R_1(n/c)\right]  =
\]
\[
\exp \left\{-n\left(\frac{\lambda}{c}-\log \frac{\lambda}{c}-1 \right)+o(n) \right\}
\exp \left\{-n\left(\frac{1}{c}-\log \frac{1}{c}-1 \right) +o(n) \right\}.
\]
Let 
\[
g_{2}(c)=\left(\frac{\lambda}{c}-\log \frac{\lambda}{c}-1\right)+ \left(\frac{1}{c}-\log \frac{1}{c}-1 \right)-\log d.
\]
Then, for all $c \in (1, \lambda)$,
 \[
\E V_n(n/c) = 
  \exp \left\{-n \cdot g_{2}(c)+o(n) \right\}.
\]
\indent Case 3: For every  $c \in [\lambda, \infty)$,
\[
u_n(n/c)
= 
\exp \left\{-n\left(\frac{1}{c}-\log \frac{1}{c}-1 \right) +o(n) \right\}.
\]
Thus, 
 \[
\E  V_n(n/c)  =
 \exp \left\{-n \cdot g_{3}(c)+o(n) \right\}
\]
where $g_{3}(c)=  \left(\frac{1}{c}-\log \frac{1}{c}-1 \right)-\log d =f(c)$. \\
\indent Define function $g(c)$ on $(0,\infty)$ by combining $g_1(c), g_2(c)$ and $g_3(c)$
on their domains.
Note that for every $c>0$ and non-negative integers $m$ and $n$
\[
u_m(m/c)u_n(n/c) \le u_{m+n}((m+n)/c).
\] 
Hence
\begin{equation} \label{TwoRich:ExpectedNumberUpperBound}
\E  V_n(n/c)\le
 \exp \left\{-n \cdot g(c) \right\}.
\end{equation}
The  function  $g(c)$ is continuosly differentiable, strictly decreasing on \\
$(0, (\lambda+1)/2 )$ and strictly increasing on $( (\lambda+1)/ 2, \infty)$ 
with the unique minimum at $c_{0}(\lambda)=\frac{\lambda+1}{2}$.
Furthermore, $g_{\lambda}(c_{0})=\log \frac{(\lambda+1)^2}{4\lambda d}$. 
It follows that  
\[
g(c_{0})>0 \mbox{  if }  \lambda> (2d-1)+\sqrt{(2d-1)^2-1},
 \]
\[ 
g(c_{0})<0 \mbox{  if }  1<\lambda< (2d-1)+\sqrt{(2d-1)^2-1}.
\]
\begin{proposition} \label{TwoIndependentRich}
Fix  $\lambda> \lambda_c=(2d-1)+\sqrt{(2d-1)^2-1}$. Almost surely, for  all sufficiently large $t$, 
\[
 R_1(t) \subset R_{\lambda}(t).
\]
\end{proposition} 
\begin{proof}
According to Proposition \ref{Rich:TreeShape}, for any $\epsilon_1 >0$ and all large $t$,
\[
\ R_{\lambda}(t) \supset D_{t(a\lambda-\epsilon_1)}.
\]
Thus, we are only interested to see what happens in 
 the region  $D^c_{t(a\lambda-\epsilon_1)}$.
Note that by (\ref{TwoRich:ExpectedNumberUpperBound}), for every $\epsilon>0$ and large $n$,\begin{equation*} \label{TwoRich:ExpectedNumberUpperBoundTimeShift}
\E \left[ \#  \left\{ x_n: x_n \in R_1(n/c) \setminus R_\lambda(n/c-1) \right\} \right]
\le
 \exp \left\{-n \cdot g(c)+n\epsilon \right\}.
\end{equation*}
For large  $t>0$ and integers $i \ge 1$, let  $c_i=i/t$.  
By Markov's inequality,
\[
\Prob \left[  (R_1(t) \setminus R_{\lambda}(t-1)) \cap D^c_{t(a\lambda-\epsilon_1)} 
\not= \emptyset
\right] \le 
\sum_{i=  [t(a\lambda -\epsilon_1)] }^{\infty} \exp \left\{ -i\cdot g(c_{i})+i\epsilon \right\} \le 
\]
\[
\le \sum_{i=  [t(a\lambda -\epsilon_1)] }^{\infty} \exp \left\{ -i g(c_{0})+i\epsilon \right\} \le C \exp \left\{ -  t(a\lambda -\epsilon_1)\cdot (g(c_{0})-\epsilon) \right\}. 
\]
A routine application of Borel-Cantelli lemma  for integer values of $t$ implies the result.
 \end{proof}
\section{The Escape model.} \label{section:Escape}
\begin{proof} [Proof of Theorem \ref{Escape:NumberOccupied} (sketch)]
Since $A(t)$ can not grow faster than the infected set in the  Richardson model with rate $1$, from the results of Section~\ref{Two independent Richardson models.}
(use (\ref{TwoRich:ExpectedNumberUpperBound}))  it follows
that for  any initial configuration $(A(0),B(0))$ and all sufficiently large $n$ 
 \[
\E  M_n(n/c)\le \exp \left\{-n \cdot g(c)+ n\epsilon \right\}.
\]
Thus, using the same lines of argument  as in Proposition  \ref{TwoIndependentRich},
it may be shown that for  all large $t$
\[
 A(t) \subset D^c_{t(r_1-\epsilon)}, \ 
\]
\[
 A(t) \subset D_{t(r_2+\epsilon)}.
\]
Furthermore, for every $c \in (r_1,r_2)$ and  any $\epsilon >0$ 
\[ \limsup \frac{1}{n}\log\left(M_{n}(n/c)\right) \le -g(c)+\epsilon. \] 
(Similarly to the proof of Proposition  \ref{Rich:NumberOccupied}, apply Markov's  inequality to get
\begin{equation*}
\Prob\left[\frac{1}{n} \log \left(M_{n}(n/c)\right) >  -g(c)+\epsilon \right] 
= \Prob\left[ M_{n}(n/c) > \exp\left\{ n(-g(c)+\epsilon) \right\} \right] \le 
\end{equation*}
\begin{equation*}
\le \E\left[ M_{n}(n/c) \right] \cdot \exp\left\{ ng(c)-n\epsilon \right\} \le \exp\left\{ -n\epsilon \right\}\end{equation*}
and  the claim follows by Borel-Cantelli lemma.)
Next, we show that on the event of survival  
\begin{equation} \label{LowerBound:A(t)}
 \liminf \frac{1}{n}\log\left(M_{n}(n/c)\right) \ge -g(c)-\epsilon.
 \end{equation} 
Note that for every  non-trivial  configuration $(A(0), B(0))$,  for all  large $t$  and  all $x \in A(t)$  we have $\mathcal{T}_{+}(x) \cap B(t)= \emptyset$.  
Furthermore, almost surely  on the event of type $1$ survival, for every integer $m$  there exist $0<t < \infty$ and a vertex $x$ with $m$-predecessor $y$  such that  
$x \in A(t)$ and  $B(t) \cap \mathcal{T}_{+}(y) = \emptyset$.
Define  $\mathcal{Z}_1(x)$ to be  a  subset of 
$C_{|x|+m} \cap \mathcal{T}_{+}(x)$ such that 
$z_1 \in \mathcal{Z}_1(x)$ if and only if  there is a {\it type $1$ path} from $(x,t)$  to $(z_1,t+m/c)$ and   there are no  directed paths from $(y,t)$ to $(x,t+m/c)$.
In general, for $k \ge 2$, $\mathcal{Z}_k(x)$ is defined as a subset of vertices in $C_{|x|+mk}  \cap \mathcal{T}_{+}(x)$ 
such that 
$z_k \in \mathcal{Z}_k(x)$ if and only if 
\begin{enumerate}
\item  $z_{k-1} \in \mathcal{Z}_{k-1}(x)$  where $z_{k-1}$ is the $m$-predecessor  of  $z_k$; 
\item there is a { \it type $1$ path} from $(z_{k-1},t+ (k-1)m/c)$ to $(z_k,t+km/c)$; 
\item there are  no directed paths from $(z_{k-2},t+ (k-1)m/c)$  to $(z_{k-1},t+km/c)$  where $z_{k-2}$ is the $m$-predecessor  of  $z_{k-1}$ .
\end{enumerate}
Let  $Z_k(x)$  be  the cardinality of the set $\mathcal{Z}_k(x)$.
It is clear from the definition that $(Z_k(x))_{k \ge 0}$ is a Galton-Watson 
process with the mean offspring number 
\[
\E[Z_1(x)]=\exp\{-mg(c)+o(m) \}.
\]
Thus, for all sufficiently large $m$,  
$\E [Z_1(x)] >1$.
 Note that $M_{|x|+mk}(t+km/c)$  dominates $Z_k(x)$,  and hence, on the event of nonextinction of  $(Z_k(x))_{k \ge 0}$,  (\ref{LowerBound:A(t)}) holds.
Observe that for every $m$,  almost surely on the event of survival of the type $1$, there are infinitely many vertices $x$ at which the Galton-Watson processes $(Z_k(x))_{k \ge 0}$ can be originated. Hence, (\ref{LowerBound:A(t)})  holds almost surely on the event of nonextinction of the first type. This finishes the proof of Theorem \ref{Escape:NumberOccupied}.
\end{proof}
\begin{proof} [Proof of Theorem \ref{Escape:MainTheorem} (sketch)]
Fix  $\lambda \in (1, \lambda_c)$, and consider any non-trivial initial configuration.
With positive probability, there exists a vertex $x$ and a Galton-Watson process
 $(Z_k(x))_{k \ge 0}$  (constructed   in the proof of Theorem  \ref{Escape:NumberOccupied})
with $\E Z_1(x) >1$. 
Hence, the Galton-Watson process survives with positive probability and 
so does type $1$. \\
\indent Consider the case   $\lambda \in (\lambda_c, \infty)$.
Since $A(t)$ is dominated  by  the  Richardson model with rate $1$, 
by  Proposition  \ref{TwoIndependentRich} in  Section 
\ref{Two independent Richardson models.}  type $1$ dies out almost surely.
\end{proof}
\section*{Acknowledgements}  
This  paper arose out of the author's  dissertation work 
completed under supervision of Steven Lalley.

\end{document}